\documentclass[10pt]{article} \usepackage{amsfonts} \usepackage{amssymb}\usepackage{amsmath}\usepackage{amsfonts}
\usepackage[russian]{babel}
\usepackage[cp1251]{inputenc}

\makeatletter
\renewcommand{\@makefnmark}{}
\makeatother
\begin{document}
\baselineskip=10pt
\pagestyle{plain}
{\Large
\makeatletter
\renewcommand{\@makefnmark}{}
\makeatother

\footnote{

\hspace{2mm}Keywords: Dirac operator, spectrum, periodic boundary conditions, inverse problem,}

\newcommand{\om}{O(\frac{1}{\mu})}
\newcommand{\lp}{L_2(0,\pi)}
\newcommand{\lpp}{L_2(\Omega)}

\newcommand{\elx}{e^{i\lambda x}}
\newcommand{\elxx}{e^{-i\lambda x}}

\newcommand{\elp}{e^{i\pi\lambda }}
\newcommand{\elpp}{e^{-i\pi\lambda }}
\newcommand{\No}{\textnumero}

\newcommand{\pp}{\prod_{n=-\infty\atop n\ne0}^\infty}
\newcommand{\sss}{\sum_{n=-\infty}^\infty}
\newcommand{\sN}{\sum_{n=-N}^N}
\newcommand{\sNN}{\sum_{|n|>N}}
\newcommand{\pps}{\sum_{n=-\infty\atop n\ne0}^\infty}
\newcommand{\sNNn}{\sum_{|n|\le N}}
\newcommand{\ppo}{\prod_{n=-\infty}^\infty}
\newcommand{\ppj}{\prod_{n=-\infty\atop n\ne0}^\infty\prod_{j=1}^2}
\newcommand{\ppjj}{\prod_{n=-\infty}^\infty\prod_{j=1}^2}

\newcommand{\ssj}{\sum_{n=-\infty\atop n\ne0}^\infty\sum_{j=1}^2}

\newcommand{\ssjj}{\sum_{n=-\infty}^\infty\sum_{j=1}^2}

\medskip
\medskip
\medskip

\centerline {\bf Characterization of the periodic and antiperiodic spectra}
 \centerline {\bf of non-self-adjoint Dirac operators}

\medskip
\medskip
\medskip
\medskip
\medskip

{\normalsize
\centerline { Alexander Makin}
\medskip
\medskip
\medskip

\medskip
\medskip
\begin{quote}{\normalsize

  ABSTRACT. The necessary and sufficient conditions are given for a sequence of complex numbers to be the periodic (or antiperiodic) spectrum of non-self-adjoint Dirac operator.}
   \end{quote}

\medskip
\medskip

\medskip
\medskip

\medskip
\medskip

\centerline {\bf 1. Intoduction}

\medskip
\medskip

 One of the important classes of inverse spectral problems is the problem of recovering a system of differential equations from spectral data.
 The solution of such problems are considered in many papers [12,18, 29-35, and the references therein]. The most studied are such problems for Dirac and Dirac type differential operator. In particular,  such problems for canonical Dirac system on a finite interval

$$
B\mathbf{y}'+V\mathbf{y} =\lambda\mathbf{y},\eqno(1.1)
 $$
where $\mathbf{y}={\rm col}(y_1(x),y_2(x))$,

 $$
 B=\begin{pmatrix}
 0&1\\
 -1&0
 \end{pmatrix},\quad  V(x)=\begin{pmatrix}
 p(x)&q(x)\\
 q(x)&-p(x)
 \end{pmatrix},
$$
in selfadjoint case have been studied in detail. In the cases of the Dirichlet and the
Newmann boundary conditions reconstruction of a continuous potential from two spectra was carried out in [6], from one spectrum and the norming constants in [5], and from the spectral function in [15]. The analogous results for Dirac operator with summable potentials were established in [1]. The case of more general separated boundary conditions was considered in [3]. In the case of unseparated boundary conditions (including periodic, antiperiodic and quasiperiodic  conditions) the indicated problem was solved in [17, 19-22].
In non-self-adjoint case
the problem of reconstructing the potential $V(x)$ from spectral data is much more complicated, since many methods
successfully used to study selfadjoint operators are inapplicable. For example, the characterization of the spectra of the periodic (antiperiodic) problem for operator (1.1) with real coefficients is given in [17] in terms of special conformal mappings, which do not exist for complex-valued potentials. The property that the eigenvalues of corresponding Dirichlet problem and Neumann problem are interlaced, which is often used to prove the solvability of the basic equation, loses its meaning in the complex case. Non-self-adjoint inverse problems for system (1.1) with different types of boundary conditions with sufficiently smooth coefficients, which, however, could have singularities were investigated in [2,8,24,28].

 Questions of uniqueness in inverse problems for operators of type (1.1) on a finite interval were studied in a lot of papers.
In particular, uniqueness of the inverse problem for general Dirac-type systems
of order 2n was established in [13, 14].
Also inverse theory was intensively developed for Dirac-type operators on the axis and semiaxis by many authors. New inverse approach to such differential  systems on the semiaxis  based on A-function concept was recently considered in [7].

 The aim of this paper is to find necessary and sufficient conditions of solvability of the periodic (antiperiodic) inverse spectral problem for system (1.1) with a nonsmooth  complex-valued potential $V(x)$.

The paper is organized as follows. Section 2 contains some basic facts and definitions related to the considered problems. In section 3 by using a modified version of the Gelfand-Levitan-Marchenko method we prove solvability of the basic equation  and establish necessary and sufficient conditions for an entire function to be the characteristic determinant of the considered problem. Further, we obtain necessary and sufficient conditions for a set of complex numbers to be the spectrum of the mentioned problem.

\medskip
\medskip

\centerline {\bf 2. Preliminaries}
\medskip
\medskip

In the present paper, we consider system  (1.1), where  complex-valued functions $p, q\in L_2(0,\pi)$ $(V\in L_2)$ with periodic (antiperiodic) boundary conditions
$$
y_1(0)-(-1)^\theta y_1(\pi)=0,\quad y_2(0)-(-1)^\theta y_2(\pi)=0,\eqno(2.1)
$$
where $\theta=0,1$.
In what follows, we introduce the Euclidean norm $\|f\|=(|f_1|^2+|f_2|^2)^{1/2}$ for vectors  $f={\rm col}(f_1,f_2)\in\mathbb{C}^2$, and set
$\langle f,g\rangle=f_1g_1+f_2g_2$. If $W$ is $2\times2$-matrix, then we set $\|W\|= \sup_{\|f\|=1}\|Wf\|$ and denote by $L_{2,2}(a,b)$ and $L_{2,2}^{2,2}(a,b)$, respectively, the spaces of 2-coordinate vector functions $f(t)={\rm col}(f_1(t),f_2(t))$
  and $2\times2$ matrix functions $W(t)$ with finite norms
$$
\|f\|_{L_{2,2}(a,b)}=\left(\int_a^b\|f(t)\|^2dt\right)^{1/2}, \quad \|W\|_{L_{2,2}^{2,2}(a,b)}=\left(\int_a^b\|W(t)\|^2dt\right)^{1/2}.
$$

The operator $\mathbb{L}\mathbf{y}=B\mathbf{y}'+V\mathbf{y}$ is regarded as a linear operator in the space $L_{2,2}(0,\pi)$
with the domain $D(\mathbb{L})=\{\mathbf{y}\in W_1^1[0,\pi]\times W_1^1[0,\pi]:\, \mathbb{L}\mathbf{y}\in L_{2,2}(0,\pi)$, $U_j(\mathbf{y})=0$ $(j=1,2)\}$.

Denote by
$$
E(x,\lambda)=\begin{pmatrix}
c_1(x,\lambda)&-s_2(x,\lambda)\\
s_1(x,\lambda)&c_2(x,\lambda)
\end{pmatrix}\eqno(2.2)
$$
the matrix of the fundamental solution system  to  equation (1.1)
with boundary condition
$
E(0,\lambda)=I
$, where
$I$ is the unit matrix,
and by
$
E_0(x,\lambda)
$
the fundamental solution system  to the equation
$B\mathbf{y}'=\lambda\mathbf{y}$
with boundary condition
$E_0(0,\lambda)=I$.
Obviously,
$$
E_0(x,\lambda)=\begin{pmatrix}
\cos\lambda x&-\sin\lambda x\\
\sin\lambda x&\cos\lambda x
\end{pmatrix}.
$$
Denote the second column of the matrix $E_0(x,\lambda)$
$$
Y_0(x,\lambda)=
\begin{pmatrix}
-\sin\lambda x\\
\cos\lambda x
\end{pmatrix}.
$$
It is well known that  the entries of the matrix
$E(x,\lambda)$ are related by the identity
$$
c_1(x,\lambda)c_2(x,\lambda)+s_1(x,\lambda)s_2(x,\lambda)=1,\eqno(2.3)
$$
which is valid for any $x, \lambda$.
The matrix $E(\pi,\lambda)$ is called the monodromy matrix of operator $\mathbb{L}\mathbf{y}$. For its entries we introduce the notation
$c_j(\lambda)=c_j(\pi,\lambda)$, $s_j(\lambda)=s_j(\pi,\lambda)$, $j=1,2$. We denote also the class of entire functions $f(z)$ of exponential type $\le\sigma$ such that $\|f\|_{L_2(R)}<\infty$ by $PW_\sigma$ .
It is known [26] that the functions $c_j(\lambda),s_j(\lambda)$ admit the representation
$$
c_j(\lambda)= \cos\pi\lambda+g_j(\lambda), \quad s_j(\lambda)=\sin\pi\lambda+h_j(\lambda),\eqno(2.4)
$$
where $g_j,h_j\in PW_\pi$, $j=1,2$. For functions of type (2.4) the following statement is true:

{\bf Lemma 1} [17] .
{\it Functions $u(\lambda)$ и $v(\lambda)$ admit the representations
$$
u(\lambda)=\sin\pi \lambda+h(\lambda),\quad v(\lambda)=\cos\pi\lambda+g(\lambda),
$$
where $h,g\in PW_\pi$,
if and only if
$$
u(\lambda)=-\pi(\lambda_0 -\lambda)\pp\frac{\lambda_n-\lambda}{n},
$$
where
$\lambda_n=n+\epsilon_n,\{\epsilon_n\}\in l_2$,
$$
v(\lambda)=\ppo\frac{\lambda_n-\lambda}{n-1/2},
$$
where
$\lambda_n=n-1/2+\kappa_n,\{\kappa_n\}\in l_2$.}

It is well known that the characteristic determinant of problem (1.1), (2.1) can be reduced to the form
$$
\Delta(\lambda)=(-1)^{\theta+1}+\frac{c_1(\lambda)+c_2(\lambda)}{2},\eqno (2.5)
$$
аnd the eigenvalues are specified by the asymptotic formulas
$$
\lambda_{n,j}=2n+\theta+\varepsilon_{n,j},\eqno (2.6)
$$
where $\{\varepsilon_{n,j}\}\in l_2$,
$n\in\mathbb{Z}$, $j=1,2$. Further $\Gamma(z,r)$ denotes a disk of radius $r$ centered at the point $z$.

Next, we establish the necessary and sufficient conditions that an entire function must satisfy in order to be the characteristic determinant of some problem (1.1), (2.1). Then, we give
an intrinsic description of sequences which are  periodic or antiperiodic spectra of operator (1.1), (2.1).

\medskip
\medskip

\centerline {\bf 3. Main Results}

\medskip
\medskip

3.1. Characteristic determinant

\medskip
\medskip

{\bf Theorem 3.1.} {\it
For a function $U(\lambda)$ to be the characteristic determinant of problem (1.1), (2.1), it is necessary and sufficient that it can be represented in the form
$$
U(\lambda)=(-1)^{\theta+1}+\cos\pi\lambda+f(\lambda),
$$
where $f\in PW_\pi$, and
$$\sss|f(n)|<\infty.\eqno (3.1)
$$}

{\bf Proof.} Necessity. Evidently, relations
(2.4), (2.5) imply that $f\in PW_\pi$. To check inequality (3.1) we consider the  monodromy matrix of problem (1.1), (2.1).
Let the corresponding function $s_2(\lambda)$ have the roots $\lambda_n$, hence by [27, Lemma 2.2],
$$
\lambda_n=n+\delta_n,\eqno (3.2)
$$
 where $\{\delta_n\}\in l_2$, $n\in\mathbb{Z}$ . Since
$$
c_j(\lambda_n)=\cos\pi\lambda_n+g_j(\lambda_n),\eqno(3.3)
$$
it follows from (2.4) and
  [27, Lemma 2.1]
that
$$
\sss|g_j(\lambda_n)|^2<\infty.\eqno(3.4)
$$
Denote
$$
\chi(\lambda)=U(\lambda)-(-1)^{\theta+1}=\cos\pi\lambda+f(\lambda).\eqno(3.5)
$$
By virtue of (2.5),
$$
c_1(\lambda_n)+c_2(\lambda_n)=2\chi(\lambda_n).
$$
It follows from (2.3) that
$
c_1(\lambda_n)c_2(\lambda_n)=1,
$
consequently the numbers $c_1(\lambda_n),c_2(\lambda_n)$ are the roots of the quadratic equation
$$
w^2-2\chi(\lambda_n)w+1=0.\eqno(3.6)
$$
Therefore we have
$$
\begin{array}{c}
c_1(\lambda_n),c_2(\lambda_n)=\chi(\lambda_n)\pm\sqrt{\chi^2(\lambda_n)-1}\\
=\cos\pi\lambda_n+f(\lambda_n)\pm\sqrt{(\cos\pi\lambda_n+f(\lambda_n))^2-1}\\
=\cos\pi\lambda_n+f(\lambda_n)\pm\sqrt{\cos^2\pi\lambda_n+2\cos\pi\lambda_nf(\lambda_n)+f^2(\lambda_n)-1}\\=
\cos\pi\lambda_n+f(\lambda_n)\pm\sqrt{2\cos\pi\lambda_nf(\lambda_n)+f^2(\lambda_n)-\sin^2\pi\lambda_n}.
\end{array}\eqno(3.7)
$$
It follows from (3.3) and (3.7) that
$$
(g_1(\lambda_n)-f(\lambda_n))^2=2\cos\pi\lambda_nf(\lambda_n)+f^2(\lambda_n)-\sin^2\pi\lambda_n,
$$
hence,
$$
2\cos\pi\lambda_nf(\lambda_n)=g_1^2(\lambda_n)-2g_1(\lambda_n)f(\lambda_n)+\sin^2\pi\delta_n. \eqno(3.8)
$$
It follows from (3.2) that for all sufficiently large $|n|$ the inequality $|\cos\pi\lambda_n|>1/2$ holds. This, together with (3.2), (3.4), and [27, Lemma 2.1] implies
$$
\sss|f(\lambda_n)|<\infty.\eqno(3.9)
$$
Since $f'\in PW_\pi$, then
$$
|f(n)|\le|f(\lambda_n)|+|f(n)-f(\lambda_n)| \le|f(\lambda_n)|+|\delta_n||\tau_n| \le|f(\lambda_n)|+(|\delta_n|^2+|\tau_n|^2)/2,
$$
where
$$
\tau_n=\max_{\lambda\in\Gamma(n,|\delta_n|)}|f'(\lambda)|.
$$
By [27, Lemma 2.1], $\{\tau_n\}\in l_2$. This and (3.9) imply (3.1).

Sufficiency. Let $f\in PW_\pi$ satisfy condition (3.1).
It follows from  the Paley-Wiener theorem and [16, Lemma 1.3.1] that
$$
\lim_{|\lambda|\to\infty}e^{-|\pi {\rm Im}\lambda|}f(\lambda)=0, \eqno(3.10)
$$
hence there exists a positive integer $N_0$ large enough that  $|f(\lambda)|<1/100$ if  ${\rm Im}\lambda=0$, $|{\rm Re}\lambda|\ge N_0$.
Let $\lambda_n$ $(n\in\mathbb{Z})$ be a strictly monotone increasing sequence of real numbers such that
for any $n\ne0$
$\lambda_n=\lambda_{-n}$,
  $|\lambda_n-(N_0+1/2)|<1/100$  if $0\le n\le N_0$, and
$\lambda_n=n$ if $n>N_0$.
Denote
$$
s(\lambda)=-\pi(\lambda_0-\lambda)\pp\frac{\lambda_n-\lambda}{n}.\eqno(3.11)
$$
It follows from Lemma 1 that
$$
s(\lambda)=\sin\pi \lambda+h(\lambda),\eqno(3.12)
$$
where $h\in PW_\pi$,
hence,
$$
|s(\lambda)|\ge C_1e^{\pi|{\rm Im}\lambda|}\eqno(3.13)
$$
if $|{\rm Im}\lambda|\ge M$,
where $M$ is sufficiently large.
It follows from (3.11) that
$$
\dot s(\lambda_0)=\pi\pp\frac{\lambda_n}{n}>0.
$$
One can readily see that the inequality $\dot s(\lambda_n)\dot s(\lambda_{n+1})<0$ holds for all $n\in\mathbb{Z}$.
It follows from two last inequalities that
$$
(-1)^n\dot s(\lambda_n)>0.\eqno(3.14)
$$
Relation (3.12) and  [27, Lemma 2.1]  imply that
$$
\dot s(\lambda_n)=\pi(-1)^n+\tau_n,\eqno(3.15)
$$
where
$\{\tau_n\}\in l_2$,
hence,
$$
\frac{1}{\dot s(\lambda_n)}=\frac{(-1)^n}{\pi}+\sigma_n,\eqno(3.16)
$$
where $\{\sigma_n\}\in l_2$.

Equation (3.6)
has the roots
$$
\begin{array}{c}
c_n^\pm=\chi(\lambda_n)\pm\sqrt{\chi^2(\lambda_n)-1}=\cos\pi\lambda_n+f(\lambda_n)\pm\sqrt{(\cos\pi\lambda_n+f(\lambda_n))^2-1}=\\
=\cos\pi\lambda_n+f(\lambda_n)\pm\sqrt{\cos^2\pi\lambda_n+2\cos\pi\lambda_nf(\lambda_n)+f^2(\lambda_n)-1}=\\=
\cos\pi\lambda_n+f(\lambda_n)\pm\sqrt{2\cos\pi\lambda_nf(\lambda_n)+f^2(\lambda_n)-\sin^2\pi\lambda_n}.
\end{array}\eqno(3.17)
$$
It follows from (3.17) that if $0<|n|\le N$ the numbers $c_n^+$
are contained within the disk $\Gamma(i,1/10)$, the numbers $c_n^-$ are contained within the disk $\Gamma(-i,1/10)$, and
if $|n|>N$ the numbers $c_n^\pm$ are contained within the disk $\Gamma(1,1/10)$ for even $n$, the numbers $c_n^\pm$ are contained within the disk $\Gamma(-1,1/10)$ for odd $n$. Denote $c_n=c_n^+$ for even $n$ and $c_n=c_n^-$ for odd  $n$. Denote also
 $$
 z_n=\frac{c_n}{\dot s(\lambda_n)}.
 $$
 It follows from (2.20) that the numbers
$z_n$  lie strictly above the line $l:{\rm Im}\lambda=-{\rm Re}\lambda$.

Evidently,
$$
\lambda_n=n+\rho_n,\eqno(3.18)
$$
where $\{\rho_n\}\in l_2$.
It follows from (3.17) and (3.18) that
$$
c_n=(-1)^n+\vartheta_n,\eqno(3.19)
$$
where $\{\vartheta_n\}\in l_2$.
Let $\beta_n=c_n-\cos\pi\lambda_n$, then $\{\beta_n\}\in l_2$. Let us consider the function
$$
g(\lambda)=s(\lambda)\sss\frac{\beta_n}{\dot s(\lambda_n)(\lambda-\lambda_n)}.
$$
By [9, p. 120] the function
$g\in PW_\pi$ and $g(\lambda_n)=\beta_n$. Denote $c(\lambda)=\cos\pi\lambda+g(\lambda)$, then $c(\lambda_n)=c_n\ne0$, hence,
the functions $s(\lambda)$ and $c(\lambda)$ have disjoint zero sets.

Denote
$$
F(x,t)=\sss\left(\frac{c_n}{\dot s(\lambda_n)}(Y_0(x,\lambda_n) Y_0^T(t,\lambda_n))-\frac{1}{\pi}Y_0(x,n) Y_0^T(t,n)\right).
$$
It follows  from [26] that
$$
\|F(\cdot,x)\|_{L_{2,2}^{2,2}(0,\pi)}+\|F(x,\cdot)\|_{L_{2,2}^{2,2}(0,\pi)}<C_2,
$$
where $C_2$ not depending on $x$.

{\bf Proposition 1.} {\it
For every $x\in[0,\pi]$ the  homogeneous equation
$$
f^T(t)+\int_0^xf^T(s)F(s,t)ds=0,\eqno(3.19)
$$
where $f(t)={\rm col}(f_1(t),f_2(t))$, $f\in L_{2,2}(0,x)$, $f(t)=0$ if $x<t\le\pi$
 has the trivial solution only.}

 Multiplying   equation (3.19)
by $\overline{f^T(t)}$ and integrating the resulting equation over segment $[0,x]$, we obtain
$$
\|f\|_{L_{2,2}(0,x)}^2+\int_0^x\langle\int_0^xf^T(s)F(s,t)ds,f^T(t)\rangle dt=0.
$$
Simple computations show
$$
\begin{array}{c}
f^T(s)F(s,t)=\sss\{z_n[f_1(s)\sin\lambda_ns\sin\lambda_nt-f_2(s)\cos\lambda_ns\sin\lambda_nt,\\ -f_1(s)\sin\lambda_ns\cos\lambda_nt
+f_2(s)\cos\lambda_ns\cos\lambda_nt]\\
-\frac{1}{\pi}[f_1(s)\sin ns\sin nt-f_2(s)\cos ns\sin nt,\\ -f_1(s)\sin ns\cos nt
+f_2(s)\cos ns\cos nt]\}\\
=\sss\{z_n[f_1(s)\sin\lambda_ns\sin\lambda_nt-f_2(s)\cos\lambda_ns\sin\lambda_nt]\\-\frac{1}{\pi}[f_1(s)\sin ns\sin nt-f_2(s)\cos ns\sin nt] ,\\ z_n[ -f_1(s)\sin\lambda_ns\cos\lambda_nt
+f_2(s)\cos\lambda_ns\cos\lambda_nt]\\
-\frac{1}{\pi}[-f_1(s)\sin ns\cos nt
+f_2(s)\cos ns\cos nt]\},
\end{array}
$$
therefore, transforming the iterated integrals into products of integrals and using the reality of all numbers $\lambda_n$, we obtain
$$
\begin{array}{c}
\int_0^x\langle f^T(s)F(s,t)ds,f^T(t)\rangle dt\\=
\sss\int_0^x\bigr(\int_0^x\{z_n[f_1(s)\sin\lambda_ns\sin\lambda_nt-f_2(s)\cos\lambda_ns\sin\lambda_nt]\\-\frac{1}{\pi}[f_1(s)\sin ns\sin nt-f_2(s)\cos ns\sin nt]\}ds\bigr)\overline{f_1(t)}dt\\+
\sss\int_0^x\bigr(\int_0^x\{z_n[-f_1(s)\sin\lambda_ns\cos\lambda_nt
+f_2(s)\cos\lambda_ns\cos\lambda_nt]\\
-\frac{1}{\pi}[-f_1(s)\sin ns\cos nt
+f_2(s)\cos ns\cos nt]\}ds\bigr)\overline{f_2(t)}dt\\=

\sss\bigr(z_n\int_0^x[f_1(s)\sin\lambda_ns-f_2(s)\cos\lambda_ns]ds\int_0^x\sin\lambda_nt\overline{f_1(t)}dt\\
-\frac{1}{\pi}\int_0^x[f_1(s)
\sin ns-f_2(s)\cos ns]ds\int_0^x\sin nt\overline{f_1(t)}dt\bigr)\\+

\sss \bigr(z_n \int_0^x[-f_1(s)\sin\lambda_ns
+f_2(s)\cos\lambda_ns]ds\int_0^x\cos\lambda_nt\overline{f_2(t)}dt\\
-\frac{1}{\pi}\int_0^x[-f_1(s)\sin ns
+f_2(s)\cos ns]ds\int_0^x\cos nt\overline{f_2(t)}dt\bigr)\\=

\sss z_n\bigr(\int_0^x[f_1(s)\sin\lambda_ns-f_2(s)\cos\lambda_ns]ds\int_0^x\sin\lambda_nt\overline{f_1(t)}dt\\+\int_0^x[-f_1(s)\sin\lambda_ns
+f_2(s)\cos\lambda_ns]ds\int_0^x\cos\lambda_nt\overline{f_2(t)}dt\bigr)\\-
\sss \frac{1}{\pi}\bigr(\int_0^x[f_1(s)\sin ns-f_2(s)\cos ns]ds\int_0^x\sin nt\overline{f_1(t)}dt\\+\int_0^x[-f_1(s)\sin ns
+f_2(s)\cos ns]ds\int_0^x\cos nt\overline{f_2(t)}dt\bigr)\\=

\sss z_n\bigr(\int_0^x[f_1(t)\sin\lambda_nt-f_2(t)\cos\lambda_nt]dt\int_0^x\sin\lambda_nt\overline{f_1(t)}dt\\+\int_0^x[-f_1(t)\sin\lambda_nt
+f_2(t)\cos\lambda_nt]dt\int_0^x\cos\lambda_nt\overline{f_2(t)}dt\bigr)\\-
\sss \frac{1}{\pi}\bigr(\int_0^x[f_1(t)\sin nt-f_2(t)\cos nt]dt\int_0^x\sin nt\overline{f_1(t)}dt\\+\int_0^x[-f_1(t)\sin nt
+f_2(t)\cos nt]dt\int_0^x\cos nt\overline{f_2(t)}dt\bigr)\\=

\sss z_n\int_0^x[f_1(t)\sin\lambda_nt-f_2(t)\cos\lambda_nt]dt\int_0^x[\overline{f_1(t)}\sin\lambda_nt-\overline{f_2(t)}\cos\lambda_nt]dt\\-
\sss\frac{1}{\pi}\int_0^x[f_1(t)\sin nt-f_2(t)\cos nt]dt\int_0^x[\overline{f_1(t)}\sin nt-\overline{f_2(t)}\cos nt]dt\\=

\sss z_n|\int_0^x\langle f(t),Y_0(t,\lambda_n)\rangle dt|^2-\sss \frac{1}{\pi}|\int_0^x\langle f(t),Y_0(t,n)\rangle dt|^2.

\end{array}
$$
It is well known that the function system $\{\frac{1}{\sqrt{\pi}}Y_0(t,n)\}$ $(n\in \mathbb{Z})$  is an orthonormal basis in $L_{2,2}(0,\pi)$, hence
it follows from the Parseval equality that
$$
\|f\|_{L_{2,2}(0,x)}^2=\sss \frac{1}{\pi}\left|\int_0^x\langle f(t),Y_0(t,n)\rangle dt\right|^2,
$$
consequently,
$$
\sss z_n\left|\int_0^x\langle f(t),Y_0(t,\lambda_n)\rangle dt\right|^2=0.
$$
Since  all the numbers $z_n$ are situated strictly in the same half-plane relative to a line which  passes through the origin, we see that
$$
\int_0^x\langle f(t),Y_0(t,\lambda_n)\rangle dt=0
$$
for all $n\in\mathbb{Z}$.
It follows from (3.12) that the function $s(\lambda)$ is a sin-type function [10], therefore [1, Lemma 5.3],  the system $Y_0(t,\lambda_n)$ is a Riesz basis of  $L_{2,2}(0,\pi)$, hence the system $Y_0(t,\lambda_n)$ is complete in
$L_{2,2}(0,\pi)$, it follows now that $f(t)\equiv0$.

By [26, Th. 5.1], the functions $c(\lambda)$  and $-s(\lambda)$ are the entries of the first line of the monodromy matrix
$$
\tilde E(\pi,\lambda)=
\begin{pmatrix}
\tilde c_1(\pi,\lambda)&-\tilde s_2(\pi,\lambda)\\
\tilde s_1(\pi,\lambda)&\tilde c_2(\pi,\lambda)
\end{pmatrix}
$$
for  problem (1.1), (2.1) with a potential $\tilde V\in L_2$, i.e.,
$$
c(\lambda)=\tilde c_1(\pi,\lambda), s(\lambda)=\tilde s_2(\pi,\lambda).\eqno(3.21)
$$
The corresponding characteristic determinant
$$
\tilde\Delta(\lambda)=(-1)^{\theta+1}+(\tilde c_1(\pi,\lambda)+\tilde c_2(\pi,\lambda))/2=(-1)^{\theta+1}+\cos\pi\lambda+\tilde f(\lambda),
$$
where $\tilde f\in PW_\pi$. It follows from  (2.3), (3.5), (3.6), (3.21) that
$$
\begin{array}{c}
\tilde\Delta(\lambda_n)=(-1)^{\theta+1}+(\tilde c_1(\pi,\lambda_n)+\tilde c_2(\pi,\lambda_n))/2\\=
(-1)^{\theta+1}+(\tilde c_1(\pi,\lambda_n)+\frac{1}{\tilde c_1(\pi,\lambda_n)})/2=(-1)^{\theta+1}+(c(\lambda_n)+\frac{1}{c(\lambda_n)})/2\\=

(-1)^{\theta+1}+\chi(\lambda_n)=U(\lambda_n).
\end{array}
$$
This implies that the function
$$
\Phi(\lambda)=\frac{U(\lambda)-\tilde\Delta(\lambda)}{s(\lambda)}=\frac{f(\lambda)-\tilde f(\lambda)}{s(\lambda)}
$$
is an entire function in the whole complex plane. Since by the Paley-Wiener theorem
$$
|f(\lambda)-\tilde f(\lambda)|<C_3e^{\pi|{\rm Im}\lambda|},\eqno(3.22)
$$
then by (3.13)
$|\Phi(\lambda)|\le C_4$ if $|{\rm Im}\lambda|\ge M$.
We denote by $H$ the set $\Gamma(N_0+1/2,1/10)\bigcup\Gamma(-N_0-1/2,1/10))\bigcup\Gamma_{|n|>N_0}(n,1/10)$.
 Since the function $s(\lambda)$ is a sin-type function [11], then $|s(\lambda)|>C_5>0$ if $\lambda\notin H$. From this inequality, (3.22) and the Maximum Principle we obtain that $|\Phi(\lambda)|<C_6$
in the strip $|{\rm Im}\lambda|\le M$, hence the function $\Phi(\lambda)$ is bounded in the whole complex plane and, by virtue of Liouville theorem, it is a constant. Let $|{\rm Im}\lambda|=M$, then it follows from (3.10) that $\lim_{|\lambda|\to\infty}(f(\lambda)-\tilde f(\lambda))=0$,
consequently $\Phi(\lambda)\equiv0$, therefore $U(\lambda)\equiv\tilde\Delta(\lambda)$.

{\bf Remark 3.1.} Necessity of condition (3.1) for the Dirac operators with skew-symmetric potentials by another method was established in [27].

3.2. Spectrum

{\bf Theorem 3.2.} {\it For a set $\Lambda$  to be the spectrum of some Dirac operator (1.1), (2.1) with a complex-valued potential $V\in L_2(0,\pi)$ it is necessary and sufficient that it consists of two sequences of eigenvalues $\lambda_{n,j}$ satisfying condition (2.6) and the inequality

$$
\sum_{k=-\infty}^\infty\left|\sss\frac{\varepsilon_{n,1}+\varepsilon_{n,2}}{2n-2k-1}\right|<\infty.\eqno(3.23)
$$

}
{\bf Proof.}
The proof of the theorem is carried out in the same lines for periodic and antiperiodic cases, and here we present the reasoning only for periodic one.

1. Sufficiency. Let two sequences $\lambda_{n,j}$ satisfy conditions
(2.6) and (3.23).
Evidently, there exists a constant $M$ such that
$$
\sup_{n,j}|\varepsilon_{n,j}|<M,\quad \ssjj|\varepsilon_{n,j}|^2<M.\eqno(3.24)
$$
It is well known that
$$
\sin\pi \lambda=\pi \lambda\pp\frac{n-\lambda}{n}=\pi \lambda\pp\left(1-\frac{\lambda}{n}\right),
$$
hence,
$$
\sin^2\frac{\pi \lambda}{2}=\frac{\pi^2\lambda^2}{4}\pp\left(\frac{2n-\lambda}{2n}\right)^2=\frac{\pi^2\lambda^2}{4}\pp\left(1-\frac{\lambda}{2n}\right)^2,
$$
therefore the function $\Delta_0(\lambda)=-1+\cos\pi\lambda$ has the representation
$$
\Delta_0(\lambda)=-\frac{\pi^2\lambda^2}{2}\pp\frac{(2n-\lambda)(2n-\lambda)}{4n^2}.\eqno(3.25)
$$
Evidently,
$$
|\Delta_0(\lambda)|< c_1e^{\pi |{\rm Im}\lambda|}.\eqno(3.26)
$$
Denote
$$
\Delta(\lambda)=-\frac{\pi^2}{2}(\lambda_{0,1}-\lambda)(\lambda_{0,2}-\lambda)\pp\frac{(\lambda_{n,1}-\lambda)(\lambda_{n,2}-\lambda)}{4n^2}.
$$

Let $f(\lambda)=\Delta(\lambda)-\Delta_0(\lambda)$. Investigation of properties of the function $f(\lambda)$ is based on the following propositions.

{\bf Proposition 2.}{\it The function $f(\lambda)$ is an entire function of exponential type not exceeding $\pi$.}

  Denote $\Gamma$ the union of the disks $\Gamma(2n,1/4)$ $(n\in Z)$. If $\lambda\notin\Gamma$, then
$$
f(\lambda)=-\Delta_0(\lambda)\left(1-\frac{\Delta(\lambda)}{\Delta_0}\right)=-\Delta_0(\lambda)(1-\phi(\lambda)),\eqno(3.27)
$$
where

$$
\begin{array}{c}
\phi(\lambda)=
\frac{(\lambda_{0,1}-\lambda)(\lambda_{0,2}-\lambda)}
{\lambda^2}\pps\frac{(\lambda_{n,1}-\lambda)(\lambda_{n,2}-\lambda)}{(2n-\lambda)(2n-\lambda)}\\
=(1-\frac{\lambda_{0,1}}{\lambda})(1-\frac{\lambda_{0,2}}{\lambda})
\ppj(1+\frac{\varepsilon_{n,j}}{2n-\lambda})=\ppjj(1+\alpha_{n,j}(\lambda)),
\end{array}
$$
where $\alpha_{0,j}(\lambda)=\frac{-\lambda_{0,j}}{\lambda}$,
 $\alpha_{n,j}(\lambda)=\frac{\varepsilon_{n,j}}{2n-\lambda}$. Let us estimate the function $\phi(\lambda)$.
It follows from (3.24) that
$$
\ssjj|\alpha_{n,j}(\lambda)|\le c_2+\ssj(|\varepsilon_{n,j}|^2+|2n-\lambda|^{-2})/2<c_3.\eqno(3.28)
$$
It is easy to see that for all $|n|>n_0$, where $n_0$ is a sufficiently large number, we have
$$
|\alpha_{n,j}(\lambda)|<1/4\eqno(3.29)
$$
for all $\lambda\notin\Gamma$.
If $|n|\le n_0$, then inequality (3.29) holds for all sufficiently large $|\lambda|$, hence
inequality (3.29) is valid for all $|\lambda|\ge C_0$.
It follows from (3.28), (3.29), and elementary inequality
$$
|\ln(1+z)|\le2|z|,\eqno(3.30)
$$
which is valid if $|z|\le1/4$ that
$$
\ssjj|\ln(1+\alpha_{n,j}(\lambda))|\le c_4.
$$
Here and throughout the following, we choose the branch of $\ln(1+z)$ that is zero for $z=0$. In view of [9, p. 433], we rewrite the last relation in the form
$$
|\phi(\lambda)|\le\ppjj|1+\alpha_{n,j}(\lambda)|\le e^{c_4}.\eqno(3.31)
$$
It follows from (3.26), (3.27), (3.31) that
$$
|f(\lambda)|<c_5e^{\pi |{\rm Im}\lambda|}\eqno(3.32)
$$
outside the domain $\Gamma'=\Gamma\cup\{|\lambda|<C_0\}$. In particular, inequality (3.32) is valid if $\lambda$ belongs lines ${\rm Im}\lambda=\pm C_0$ and vertical segments with vertexes $(2k-1, -C_0),(2k-1, C_0)$, $|2k-1|>C_0$, $k\in\mathbb{Z}$. By the Maximum Principle, inequality (3.32) holds in whole complex plane, hence the function $f(\lambda)$ is an entire function of exponential type not exceeding $\pi$.

{\bf Proposition 3.}{\it The function $f$ belongs to $PW_\pi$.}

 Denote
$$
W(\lambda)=\ln\phi(\lambda)=\ssjj\ln(1+\alpha_{n,j}(\lambda)),
$$
then
$$
f(\lambda)=-\Delta_0(\lambda)\left(1-e^{W(\lambda)}\right).\eqno(3.33)
$$
 Let us estimate the function $W(\lambda)$ if $\lambda\notin\Gamma'$.
It follows from (3.24), (3.29), (3.30) that
$$
\begin{array}{c}
|W(\lambda)|\le\ssjj|\ln(1+\alpha_{n,j}(\lambda))|\le2\ssjj|\alpha_{n,j}(\lambda)|\\\le
\frac{2M}{|\lambda)|}+
\ssj(\frac{|\varepsilon_{n,j}|^2}{10M}+\frac{10M}{|2n-\lambda|^2})\le\frac{2M}{|\lambda|}+1/10+20M\sum_{n=0}^\infty\frac{1}{n^2+
|{\rm Im}\lambda|^2}\\\le
\frac{2M}{|\lambda|}+1/10+20M(\frac{2}{|{\rm Im}\lambda|^2}+\int_1^\infty\frac{dx}{x^2+|{\rm Im}\lambda|^2})\le\frac{2M}{|{\rm Im}\lambda|}
+1/10+20M(\frac{2}{|{\rm Im}\lambda|^2}+\frac{\pi}{2|{\rm Im}\lambda|}).
\end{array}
$$
The last inequality implies that
$$
|W(\lambda)|<1/4\eqno(3.34)
$$
if $|{\rm Im}\lambda|\ge M_1=10(\pi+2+22M)+C_0$.
Then from the trivial inequality
 $$
\frac{|z|}{2}\le |1-e^z|\le2|z|,\eqno(3.35)
 $$ which holds for $|z|\le1/4$, we obtain the inequality $|1-e^{W(\lambda)}|\le2|W(\lambda)|$, which, together with (3.26) and (3.33) implies that
$$
|f(\lambda)|\le c_6|W(\lambda)|\eqno(3.36)
$$
for $\lambda\in l$, where $l$ is the line ${\rm Im}\lambda=M_1$. Let us prove that
$$
\int_l|W(\lambda)|^2d\lambda<\infty.\eqno(3.37)
$$
From the elementary inequality $|\ln(1+z)-z|\le|z|^2$ true for $|z|\le1/2$, we obtain
$$
\ln(1+z)-z=r(z),
$$
where $|r(z)|\le|z|^2$, hence,
$$
W(\lambda)=S_1(\lambda)+S_2(\lambda),\eqno(3.38)
$$
where
$$
S_1(\lambda)=\ssjj\alpha_{n,j}(\lambda),
$$

$$
|S_2(\lambda)|\le\ssjj|\alpha_{n,j}(\lambda)|^2.
$$
Evidently,
$$
|W(\lambda)|\le|S_1(\lambda)|+|S_2(\lambda)|.\eqno(3.39)
$$

Set
$$
I_m=\int_l|S_m(\lambda)|^2d\lambda
$$
$(m=1,2)$. First consider the integral $I_1$. It follows from [25] that
$$
I_1=\int_l\left|\sss\frac{\varepsilon_{n,1}+\varepsilon_{n,2}}{2n-\lambda}\right|^2d\lambda<\infty.\eqno(3.40)
$$
It is readily seen that
$$
|S_2(\lambda)|\le\sss\frac{|\varepsilon_{n,1}|^2+|\varepsilon_{n,2}|^2}{|2n-\lambda|^2}<c_7,
$$
hence,
$$
\begin{array}{c}
I_2\le c_7\int_l\left(\sss\frac{|\varepsilon_{n,1}|^2+|\varepsilon_{n,2}|^2}{|2n-\lambda|^2}\right)d\lambda=
c_8\sss(|\varepsilon_{n,1}|^2+|\varepsilon_{n,2}|^2)\int_l\frac{d\lambda}{|2n-\lambda|^2}\\<
c_9\sss(|\varepsilon_{n,1}|^2+|\varepsilon_{n,2}|^2)<c_{10}.
\end{array}\eqno(3.41)
$$
Relations (3.39-3.45) imply (3.37). It follows from (3.36), (3.37),  and [22, p.  115] that
$$
\int_R|f(\lambda)|^2d\lambda<\infty.\eqno(3.42)
$$

{\bf Proposition 4.}{\it The function $f(\lambda)$ satisfies condition (3.1).}

Obviously,
$\Delta_0(2k)=0$, hence $f(2k)=\Delta(2k)$. Since the function $\Delta(\lambda)$ is bounded in the strip $|{\rm Im}\lambda|\le1$, and for all sufficiently large $|k|$ the inequality $|\varepsilon_{k,j}|<1/2$ takes place, by the Maximum Principle, we have
$$
|f(2k)|=|\Delta(2k)|\le|\varepsilon_{k,1}||\varepsilon_{k,2}|\max_{|2k-\lambda|=1}\left|\frac{\Delta(\lambda)}{(\lambda_{k,1}-\lambda)(\lambda_{k,2}-\lambda)}\right|\le c_{11}(|\varepsilon_{k,1}|^2+|\varepsilon_{k,2}|^2).\eqno(3.43)
$$

Let us estimate $|f(2k+1)|$. Obviously, $\Delta_0(2k+1)=-2$.
Denote
 $$\epsilon_n=\max(|\varepsilon_{n,1}|,|\varepsilon_{n,2}|).
$$
 There exists  a number $n_0>0$ such that
$$
\sum_{|n|>n_0}\epsilon_n^2<1/1000,
$$
and
for any $|n|>n_0$ the inequality
$
\epsilon_n^{2/3}<1/1000
$ holds.
Let $\lambda\notin\Gamma'$. Supplementary suppose that
$$
|\lambda|>M_2=1000(2n_0+1)n_0M.
$$
Then, using the well known inequality $ab\le \frac{a^p}{p}+\frac{b^q}{q}$ $(a,b>0, p,q>1, 1/p+1/q=1)$, we obtain
$$
\begin{array}{c}
\ssjj|\alpha_{n,j}(\lambda)|\le2(\sum_{|n|\le n_0}\frac{\epsilon_n}{|2n-\lambda|}+\sum_{|n|> n_0}\frac{\epsilon_n}{|2n-\lambda|})\\
\le2M\sum_{|n|\le n_0}\frac{1}{|2n-\lambda|}+2\sum_{|n|> n_0}(\epsilon_n^2+\frac{\epsilon_n^{2/3}}{|2n-\lambda|^{4/3}})\\
\le1/50+1/500\sum_{n=1}^\infty\frac{1}{n^{4/3}}<1/10,
\end{array}\eqno(3.44)
$$
hence inequality (3.34) is valid
for all $\lambda$ belonging to the considered domain.
Arguing as above, we see that
$$
|f(\lambda)|\le c_{12}\left(|\ssjj\alpha_{n,j}(\lambda)|+\ssjj|\alpha_{n,j}(\lambda)|^2\right).
$$
The last inequality implies that for all $|2k+1|>k_0$, where $k_0=\max(C_0, M_2)$,

$$
|f(2k+1)|\le c_{13}\left(\left|\sss\frac{\varepsilon_{n,1}+\varepsilon_{n,2}}{2n-2k-1}\right|+
\sss\frac{|\varepsilon_{n,1}|^2+|\varepsilon_{n,2}|^2}{|2n-2k-1|^2}\right).\eqno(3.45)
$$
Clearly,
$$
\begin{array}{c}
\sum_{k=-\infty}^\infty\sss\frac{|\varepsilon_{n,1}|^2+|\varepsilon_{n,2}|^2}{|2n-2k-1|^2}=
\sss(|\varepsilon_{n,1}|^2+|\varepsilon_{n,2}|^2)\sum_{k=-\infty}^\infty\frac{1}{|2n-2k-1|^2}\\
<c_{14}\sss(|\varepsilon_{n,1}|^2+|\varepsilon_{n,2}|^2)<c_{15}.
\end{array}\eqno(3.46)
$$
It follows from (3.23),(3.43), (3.45), (3.46) that (3.1) holds, hence the function $f(\lambda)$ satisfies all conditions of Theorem 3.1, and the function $\Delta(\lambda)$ is the characteristic determinant of some problem (1.1), (2.1).

2. Necessity. If a set $\{\Lambda\}$ is the spectrum of a Dirac operator (1.1), (2.1), then relation (2.6) takes place [4, Th. 6.5].
 Let us prove that condition (3.23) holds.
Since $f(\lambda)=\Delta(\lambda)-\Delta_0(\lambda)$, then by Theorem 3.1  relation (3.1) is valid.

Let $\lambda=2k+1$, $k\in\mathbb{Z}$, $|2k+1|>k_0$, hence inequality (3.44) holds. Since $\Delta_0(2k+1)=-2$, it follows from (3.33) and (3.35) that
$$
|W(2k+1)|\le|f(2k+1)|.
$$
This, together with (3.38) implies
$$
|S_1(2k+1)|\le|f(2k+1)|+\ssjj|\alpha_{n,j}(2k+1)|^2.\eqno(3.47)
$$
Using (3.46), we find that
$$
\ssjj|\alpha_{n,j}(2k+1)|^2<c_{16}.\eqno(3.48)
$$
It follows from (3.47), (3.48), and (3.1) that

$$
\sum_{|2k+1|>k_0}|S_1(2k+1)|<c_{17}.
$$
It is easy to see that
$$
\sum_{|2k+1|\le k_0}|S_1(2k+1)|<k_0c_{18}.
$$
The last two inequalities imply (3.23).

3. Example. We give an example when (2.6) holds but (3.23) does not hold. Let $\varepsilon_{n,1}=\varepsilon_{n,2}=1/m$ if $n=2^m$ and
$\varepsilon_{n,1}=\varepsilon_{n,2}=0$ if $n\ne2^m$, $m=1,2,\ldots$, $n\in\mathbb{Z}$. Denote
$$
\gamma_k=\sss\frac{\varepsilon_{n,1}+\varepsilon_{n,2}}{2n-2k-1},\quad k\in\mathbb{Z},
$$
hence,
$$
\gamma_k=2\sum_{m=1}^\infty\frac{1}{m(2^{m+1}-2k-1)}.
$$
Let $k=2^p$, $p\in\mathbb{N}$, then
$$
\gamma_{2^p}=2\sum_{m=1}^\infty\frac{1}{m(2^{m+1}-2^{p+1}-1)}=-\frac{2}{p}+\sigma_{p,1}+\sigma_{p,2},
$$
where
$$
\sigma_{p,1}=
2\sum_{m=1}^{p-1}\frac{1}{m(2^{m+1}-2^{p+1}-1)},\quad\sigma_{p,2}=2\sum_{m=p+1}^\infty\frac{1}{m(2^{m+1}-2^{p+1}-1)}.
$$
A simple computation shows that
$$
|\sigma_{p,1}|\le\frac{1}{2^p}\sum_{m=1}^{p-1}\frac{1}{m(1-2^{m-p}-1/2^{p+1})}\le\frac{4}{2^p}\sum_{m=1}^{p-1}\frac{1}{m}\le\frac{4(1+\ln p)}{2^p}\eqno(3.49)
$$
and
$$
|\sigma_{p,2}|\le2\sum_{l=1}^\infty\frac{1}{(l+p)(2^{l+p+1}-2^{p+1}-1)}
\le\frac{1}{p2^p}\sum_{l=1}^\infty\frac{1}{2^l-1-1/2^{p+1}}\le\frac{4}{p2^p}.\eqno(3.50)
$$
It follows from (3.49), (3.50) that  $|\sigma_{p,1}+\sigma_{p,2}|\le1/p$ if  $p\ge10$, hence $|\gamma_{2^p}|>1/p$, therefore the series
$$
\sum_{k\in\mathbb{Z}}|\gamma_k|
$$
diverges.

\medskip
\medskip

\centerline{\bf References}

\medskip
\medskip

[1] A. Albeverio, R. Hryniv, and  Ya. Mykytyuk,  {\it Inverse spectral problems for Dirac operators with summable potentials,} Russ. J. Math. Phys. {\bf 12}   (2005), no. 4, 406-423.

[2] N. Bondarenko, S. Buterin {\it An inverse spectral problem for integro-differential Dirac operators with general convolution kernels,}
Appl. Analysis {\bf 99} (2020), no. 4, 700-716.

[3] Vasil B. Daskalov and Evgeni Kh. Khristov
{\it
Explicit formulae for the inverse problem for the regular
Dirac operator,}
Inverse Problems {\bf 16} (2000), 247–258.

[4] P. Djakov and  B. Mityagin, {\it Unconditional Convergence of Spectral Decompositions of 1D Dirac operators with Regular Boundary Conditions,}  Indiana Univ. Math. J. \textbf{61} (2012), 359-398.

[5] T.T. Dzabiev, {\it The inverse problem for the Dirac equation with a singularity,} Acad. Nauk. Azerbaidzan. SSR.  Dokl. {\bf 22}, (1966), no. 11, p. 8-12. (in Russian).

[6] M.G. Gasymov, T.T. Dzabiev, {\it Solution of the inverse problem by two spectra for the Dirac equation on a finite interval,} Acad. Nauk. Azerbaidzan. SSR.  Dokl. {\bf 22}, (1966), no. 7, p. 3-7. (in Russian).

[7] F. Gesztesy, A. Sakhnovich,  {\it The inverse approach to Dirac-type systems based on the A-function concept,} Journal of Functional Analysis, Volume 279, Issue 6, 1 October 2020, 108609.

[8] O. Gorbunov O., V. Yurko, {\it Inverse problem for Dirac system with singularities in interior points,} Anal Math. Phys. {bf 6},  (2016), no. 1, 1-29.

[9] M.A. Lavrentiev, B.V. Shabat  {\it Methods of Theory of Complex Variable,} Nauka, Moscow, 1973 (in Russian).

[10] B. Ya. Levin, {\it Lectures on Entire Functions,} Amer. Math. Soc. Transl. Math. Monographs {\bf 150}, Amer. Math. Providence, RI, 1996.

[11] B. Ya. Levin, I. V. Ostrovskii, {\it On small perturbations of the set of zeros of functions of sine type,} Math. USSR-Izv. {\bf 14}, (1980), no. 1, 79–101.

[12] B.M. Levitan, I.S. Sargsyan, {\it Sturm–Liouville and Dirac operators,} Kluwer
Academic Publishers, Dordrecht, 1991.

[13] M. M. Malamud, {On Borg-type theorems for first-order systems on a finite interval,} Funktsional.
Anal. i Prilozhen. 33 (1999), no. 1, 75-80 (Russian); Engl. transl.: Funct. Anal. Appl. 33 (1999),
no. 1, 64-68.

[14] M. M. Malamud, {\it Uniqueness questions in inverse problems for systems of differential equations on a finite interval,}
Trans. Moscow Math. Soc. {\bf 60} (1999), 173–224.

[15] S.G. Mamedov, {\it The inverse boundary problem on a finite interval for the Dirac's systems of equations,} Azerbaidzan Gos. Univ. Uchen. Zap. Ser. Fiz-Mat. Nauk.  (1975), no. 5, 61-67 (in Russian).

[16] V. A. Marchenko, {\it Sturm–Liouville operators and their applications,}
Birkhauser Verlag, Basel, 1986.

[17] T.V. Misyura, {\it Characterization of spectra of periodic and anti-periodic problems generated by Dirac's operators. II,} Theoriya functfii, funct. analiz i ikh prilozhen. {\bf 31} (1979), 102-109.

[18] Y. V. Mykytyuk,  D.V. Puyda, {\it Inverse spectral problems for Dirac operators on a finite interval,} J. Math. Anal. Appl.
{\bf 386} (2012), no. 1, 177–194.

[19] I.M. Nabiev,  {\it Solution of the class of inverse problems for the Dirac operators,} Trans. Acad. Sci. Azerb. Ser. Phys.-Tech. Math. Sci. Math. Mech. {\bf 21} (2001), 146-157.

[20] I.M. Nabiev,{\it  Characteristic of spectral data of Dirac operators,} Trans. Acad. Sci. Azerb. Ser. Phys.-Tech. Math. Sci. Math. Mech. {\bf 24} (2004), 161-166.

[21] I.M. Nabiev, {\it Solution of the quasiperiodic problem for the Dirac system,} Math. Notes, {\bf 89} (2011), no. 6, 845-852.

[22] I.M. Nabiev, {\it The inverse periodic problem for the Dirac operator,} Proceeding of IMM of NAS of Azerbaijan, {\bf XIX} (2003), 177-180.

[23] S.M. Nikolskii,  {\it Approximation of Functions of Several Variables and Embedding Theorems,} Nauka, Moscow, 1977 (in Russian).

[24] W. Ning, {\it An inverse spectral problem for a nonsymmetric differential operator: Reconstruction of eigenvalue problem,} J. Math. Anal. Appl. {\bf 327} (2007), no. 2, 1396-1419.

[25] J.-J. Sansug, V. Tkachenko, {\it Characterization of the periodic and antiperiodic spectra of nonselfadjoint Hill's operators,} Oper. Theory Adv. and Appl. {\bf 98} (1997), 216-224.

[26] V. Tkachenko, {\it Non-self-adjoint Periodic Dirac Operators,} Oper. Theory: Adv. and Appl. {\bf 123} (2001), 485-512.

[27] V. Tkachenko, {\it Non-self-adjoint Periodic Dirac Operators with finite-band spectra,} Int. Equ. Oper. Theory {\bf 36} (2000), 325-348.

[28] C.-Fu Yang and V. Yurko,  {\it Recovering Dirac operator with nonlocal boundary conditions,} J. Math. Anal. Appl. {\bf 440} (2016), no 1, 156-166.

[29]  V.A. Yurko, {\it Method of spectral mappings in the inverse problem theory, Inverse and Ill-posed
Problems Series,} VSP, Utrecht,  2002.

[30] V.A. Yurko, {\it Introduction to the Theory of Inverse Spectral Problems,} Fizmatlit, Moscow, 2007 (in Russian).

[31] V.A. Yurko  {\it Inverse spectral problems for differential operators and their applications,} Gordon and
Breach Science Publishers, Amsterdam, 2000.

[32] V.A.Yurko,  {\it An inverse spectral problem for singular non-selfadjoint differential systems,}
 Sbornik: Mathematics {\bf 195} (2004), no. 12, 1823–1854.

[33]  V.A. Yurko {\it Inverse spectral problems for differential systems on a finite interval,} Results Math. {\bf 48} (2005), no. (3–4),
371–386.

[34] V.A. Yurko {\it An inverse problem for differential systems on a finite interval in the case of multiple roots of the
characteristic polynomial,}
 Diff. Eqns. {\bf 41} (2005), no. 6, 818–823.

[35]  V.A. Yurko, {\it An inverse problem for differential systems with multiplied roots of the characteristic polynomial,} J.
Inv. Ill-Posed Probl. {\bf 13} (2005), no. 5, 503–512.

\medskip
\medskip
\medskip
\medskip

Russian Technological University, Prospect Vernadskogo 78, Moscow, 119454, Russia

\medskip
\medskip

email: alexmakin@yandex.ru

\end{document}